\documentclass[12pt]{article}
\usepackage[utf8]{inputenc}
\usepackage[english,russian]{babel}
\usepackage{graphicx}
\usepackage{vmumath}
\usepackage{amsmath}
\usepackage{amsthm}
\usepackage{avo_book}
\usepackage{tikz}
\usepackage{caption}
\usepackage{booktabs}
\usepackage{subcaption}
\usepackage[a4paper,left=15mm,right=15mm,top=30mm,bottom=20mm]{geometry}
\noaftermath
\begin{document}

\bibliographystyle{unsrt}

\title{Перечисление регулярных карт\\
на поверхностях заданного рода}

\author{
Е.С.Краско, \qquad  А.В.Омельченко\\
\small Санкт-Петербургский Академический университет РАН\\
\small Улица Хлопина, д.8, корп.3, Санкт-Петербург, 194021, Россия\\
\small\tt \{krasko.evgeniy, avo.travel\}@gmail.com
}

\begin{abstract}
В представленной работе перечисляются помеченные и непомеченные $d$-регулярные карты на двумерных ориентированных поверхностях произвольного рода $g$. Отдельно и более подробно рассматривается случай $d$-регулярных карт с одной гранью.
\end{abstract}

\maketitle

\section{Введение}

Под картой $\M$ на поверхности $X$ рода $g$ мы будем понимать вложение связного мультиграфа $G$ в  замкнутую ориентированную поверхность $X$, такое, что $X\setminus G$ представляет собой набор из $f$ связных областей, гомеоморфных открытому диску. Такие области представляют собой грани карты $\M$, имеющей $v$ вершин (точек на поверхности $X$) и $n$ ребер (несамопересекающихся кривых на поверхности, не имеющих общих точек, отличных от вершин $v$). Под $d$-регулярными картами понимаются карты, в которых каждая из вершин имеет одну и ту же степень $d$. Карта $\M$ с выделенным и ориентированным ребром называется корневой картой. 

Одной из первых работ, посвященных перечислению корневых карт на поверхности рода $g$, была работа Уолша и Лемана \cite{Walsh_Lehman}. В этой работе авторы, используя подход Татта  \cite{Tutte_Census}, \cite{Tutte_Enum}, к перечислению планарных карт, получили рекуррентные соотношения для количества корневых карт на поверхности рода $g$ и подсчитали их количество для фиксированного числа ребер. Кроме того, им удалось получить точное аналитическое выражение для количества $\epsilon_g(n)$ карт с единственной гранью, построенных на поверхности рода $g$ и имеющих $n$ ребер. В середине восьмидесятых годов прошлого века Харер и Загир с помощью метода матричных интегралов получили красивое рекуррентное соотношение для чисел $\epsilon_g(n)$ \cite{HarerZagier}. В достаточно свежей работе \cite{Chapuy} Чапуи с использованием оригинальных топологических соображений предложил изящную комбинаторную интерпретацию этих рекуррентных соотношений. Там же он показал, как эти соображения можно использовать для перечисления специальных подклассов карт с единственной гранью, в частности, для перечисления $3$-регулярных (кубических) карт. 

Параллельно с отмеченными выше исследованиями в восьмидесятые годы двадцатого века появились работы \cite{Liskovets_85}, \cite{Wormald}, посвященные перечислению так называемых непомеченных (некорневых) карт на плоскости. Идеи, лежащие в работе \cite{Liskovets_85}, были существенным образом усилены в статье Медных и Недели \cite{Mednykh_Nedela}. В этой статье для перечисления некорневых карт на ориентированных поверхностях заданного рода $g$ было введено важное понятие фактор-карт, или карт на орбифолдах. Используя это понятие, Медных и Неделе удалось вывести явную аналитическую формулу, позволяющую перечислять некорневые карты на замкнутых ориентированных поверхностях заданного рода $g$ в зависимости от числа ребер $n$. 

Представленная работа посвящена перечислению $d$-регулярных карт на поверхностях заданного рода. В первой части статьи мы, используя формулы, полученные Уолшем и Леманом, а также подход, описанный в работе \cite{Mednykh_Nedela}, предъявляем алгоритм перечисления непомеченных $d$-регулярных карт с единственной гранью на поверхности заданного рода $g$. В случае простых $d>2$ мы приводим явные аналитические формулы для подсчета таких карт, обобщающие соотношения, полученные в \cite{Bacher_Vdovina} для случая $d=3$. Во второй части данной статьи мы выводим формулы, позволяющие перечислить как корневые, так и некорневые $d$-регулярные карты на ориентированных поверхностях. В конце статьи мы приводим таблицы для числа таких карт при различных значениях $g$ и $d$.  

\section{Постановка задачи}

Формула Эйлера связывает параметры $n$, $v$ и $f$ для заданной карты $\M$ на поверхности рода $g$:
\begin{equation}
\label{eq:Euler_multicell}
v-n+f=2-2g.
\end{equation}
В случае карты с единственной гранью последнюю формулу можно упростить: 
\begin{equation}
\label{eq:Euler_unicell}
n-v=2g-1.
\end{equation} 
Кроме того, для $d$-регулярных карт из так называемой леммы о рукопожатиях (сумма степеней вершин равна удвоенному количеству ребер) мы получаем равенство вида
\begin{equation}
\label{eq:handsh_lemma}
d\cdot v=2n.
\end{equation}

Анализ формул (\ref{eq:Euler_unicell}) и (\ref{eq:handsh_lemma}) показывает, что карты с единственной гранью существуют не при всех значениях параметра $d$. Действительно, из уравнения Эйлера (\ref{eq:Euler_unicell}) следует, что в случае $f=1$ разность $n-v$ четной быть не может. Подставляя в (\ref{eq:handsh_lemma}) вместо $n$ выражение вида
$$
n=2i+v,\qquad i=0,1,2,\ldots,
$$
мы получаем, что $d$-регулярные карты с единственной гранью невозможно построить в случае, когда
$$
d=2+4i,\qquad i=0,1,2,\ldots
$$
Кроме того, только $3$ и $4$-регулярные карты можно нарисовать на поверхностях произвольного рода $g\geq 0$. Для других значений параметра $d$ система уравнений (\ref{eq:Euler_unicell}), (\ref{eq:handsh_lemma}) имеет целочисленные решения только лишь в случае, когда $4g-2$ делится на $d-2$. К примеру, в случае $d=5$ карты с единственной гранью существуют только на поверхностях рода $g=2+3i$, $i=0,1,2,\ldots$. Для карт с произвольными значениями параметра $f$ такие ограничения на род поверхности $X$ отсутствуют.

\section{Основные принципы перечисления непомеченных карт}

Основная цель данной статьи состоит в перечислении непомеченных $d$-регулярных карт на поверхностях. Сделать это непосредственно обычно довольно затруднительно вследсвтие наличия у таких карт нетривиальных групп автоморфизмов. Поэтому традиционно на первом этапе решается задача перечисления помеченных карт, а затем с помощью, например, леммы Бернсайда, подсчитывается количество непомеченных карт. 

Нам будет удобно в качестве помеченной карты $\M$ рассматривать карту, в которой выделено одно из ее полуребер (или пара ``вершина-ребро''). Такую карту мы будем называть корневой. Так как любой автоморфизм карты $\M$, оставляющий корневое полуребро на месте, тривиален \cite{Tutte_Enum}, то при перечислении корневых карт мы можем не учитывать симметрии поверхности $X$. 

После получения количества корневых $d$-регулярных карт можно попытаться, используя описанный в работе \cite{Mednykh_Nedela} подход, подсчитать количество соответствующих им непомеченных карт на поверхностях рода $g$. Основной инструмент для перечисления непомеченных объектов --- Лемма Бернсайда --- позволяет выразить это число, зная количество помеченных объектов, являющихся неподвижными точками для элементов некоторой группы $S$, орбиты действия которой мы и считаем непомеченными объектами. В \cite{Mednykh_Nedela} было показано, что при рассмотрении карт на многообразиях в качестве элементов $\sigma$ группы $S$ нам достаточно рассмотреть множество периодических сохраняющих ориентацию гомеоморфизмов поверхности, на которую уложена карта. Таким образом мы сталкиваемся со следующей задачей: по заданному гомеоморфизму описать и перечислить все карты, которые переходят в себя под его действием. Для этого мы воспользуемся идеей так называемых фактор-карт \cite{Mednykh_Nedela}.

Проиллюстрируем эту идею на примере карт $\M$ на торе. Рассмотрим представление тора в виде $2n$-угольника (например, квадрата с отождествленными противоположными сторонами). Предположим, что имеется возможность так нарисовать карту $\M$ на квадрате, представляющем тор, чтобы она переходила в себя при вращениях этого квадрата на угол $2\pi/L$. В этом случае говорят, что данная карта допускает автоморфизм периода $L$ (обладает симметрией с периодом $L$).

Заметим теперь, что при таких вращениях все множество точек квадрата разбивается на два подмножества --- бесконечное подмножество точек общего положения и конечное множество особых (критических) точек. Любая точка общего положения характеризуется тем, что мощность ее орбиты при вращении на угол $2\pi/L$ совпадает с $L$. Особая точка --- это точка, мощность орбиты которой строго меньше $L$.

\begin{figure}[ht]
\centering
	\begin{subfigure}[b]{0.45\textwidth}
	\centering
    		\includegraphics[scale=1.3]{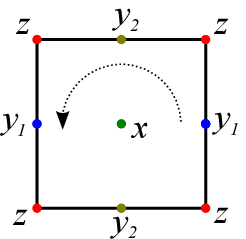}
		\caption{}
	\end{subfigure}	
	\begin{subfigure}[b]{0.45\textwidth}
	\centering
    		\includegraphics[scale=1.3]{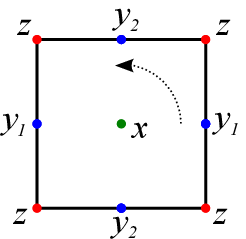}
		\caption{}
	\end{subfigure}
\caption{Автоморфизмы тора}
\label{fig:square}
\end{figure}

В качестве первого примера рассмотрим вращение квадрата, представляющего тор, на угол в $180^\circ$ ($L=2$) (рис. \ref{fig:square},a). У точек общего положения мощность орбиты под действием группы вращений $\Z_2$ совпадает с порядком этой группы и равна двум. Однако, помимо этих точек, у тора существуют четыре особые точки, остающиеся неподвижными при таких вращениях: точка $x$, отвечающая центру квадрата, точка $z$, отвечающая четырем вершинам квадрата (при склеивании противоположных сторон квадрата все эти четыре точки переходят в единственную точку $z$ на торе), а также две точки $y_1$, $y_2$, отвечающие серединам сторон квадрата. Иными словами, мощность орбит этих точек одинакова и равна единице. 

Рассмотрим теперь вращение тора на угол в $90^\circ$ ($L=4$) (рис. \ref{fig:square},b). При таком вращении у нас по-прежнему имеются те же самые четыре особых точки. Две из них --- точки $x$ и $z$ --- по-прежнему остаются неподвижными под действием таких вращений (то есть мощность соответствующих этим точкам орбит равна единице). Точки же $y_1$ и $y_2$ переходят при таких вращениях друг в друга. При этом мощность общей для этих двух точек орбиты равна двум. Отметим для полноты картины, что мощность орбиты любой точки общего положения равна четырем.

Заметим теперь, что с топологической точки зрения любой такой тор $X$ можно рассматривать как некоторое $L$-листное разветвленное накрытие сферы $O$, имеющей несколько так называемых точек ветвления (branch points). Многообразие $O$ с конечным числом точек ветвления в общем случае называется $g$-admissible orbifold ($g$-допускающий орбифолд), где $g$ --- род накрывающей $O$ поверхности $X$. Про исходное многообразие $X$ (в данном случае --- тор) говорят, что оно накрывает многообразие $O$. Само же накрытие является разветвленным в следующем смысле: все точки пространства $O$, кроме конечного числа, имеют одинаковое количество $L$ прообразов на многообразии $X$. У точек же ветвления количество прообразов строго меньше $L$. Соответствующие этим точкам прообразы и являются особыми (критическими) точками на накрывающем $O$ многообразии $X$. Количество прообразов любой точки ветвления совпадает с мощностью орбиты любого такого прообраза под действием группы $\Z_L$. Саму же точку ветвления можно характеризовать с помощью так называемого индекса ветвления (branch index) --- числа $m$, равного периоду (или количеству листов в накрытии) $L$, деленному на количество прообразов данной точки. Наконец, все разветвленное $L$-листное накрытие многообразия $O$ многообразием $X$ можно описать с помощью так называемой сигнатуры
$$
(g,\mathfrak{g},L,[m_1,\ldots,m_r]),\qquad 1<m_1\leq\ldots\leq m_r,
$$
в которой $g$ есть род накрывающей поверхности $X$, $\mathfrak{g}$ --- род орбифолда $O$, $L$ --- количество листов в накрытии, $r$ --- количество точек ветвления, $m_i$ --- индексы ветвления.

Разберем с этой точки зрения разобранные выше примеры. В случае $L=2$ мы можем рассматривать двулистное накрытие тором $X$ сферы $O$ c четырьмя точками ветвления индексов $m_i=2$. Следовательно, сигнатура такого накрытия имеет вид
$$
(1,0,2,[2,2,2,2])\equiv (1,0,2,[2^4]).
$$
В случае $L=4$ мы имеем четырехлистное накрытие сферы тором с тремя точками ветвления на сфере и с четырьмя критическими точками на торе. Сигнатура описываемого разветвленного накрытия имеет вид
$$
(1,0,4,[2,4,4])\equiv (1,0,4,[2,4^2]).
$$

Вернемся теперь к карте $\M$, нарисованной на торе и переходящей в себя при описанных выше вращениях на угол $2\pi/L$. Рассмотрим $1/L$-ю часть этой карты. Понятно, что исходную карту мы можем получить, склеивая $L$ таких частей между собой. Каждая $1/L$-я часть исходной карты $\M$ представляет собой так называемую фактор-карту (quotient map) $\mathfrak{M}$, нарисованную на сфере с точками ветвления. Фактор-карта отличается от обычной карты тем, что у нее, помимо обычных ребер, могут появляться так называемые висячие полуребра. Кроме того, наличие точек ветвления на сфере усложняет соответствие между исходной картой $\M$ и фактор-картой $\mathfrak{M}$. 

\begin{figure}[ht]
\centering
	\begin{subfigure}[b]{0.3\textwidth}
	\centering
    		\includegraphics[scale=1.3]{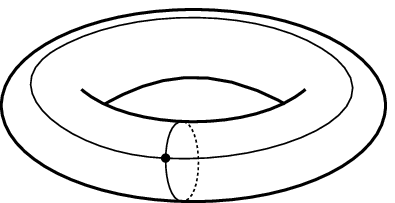}
		\caption{}
	\end{subfigure}
	\begin{subfigure}[b]{0.3\textwidth}
	\centering
    		\includegraphics[scale=1.3]{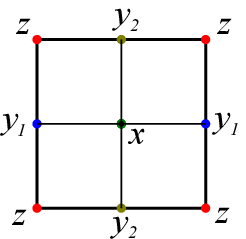}
		\caption{}
	\end{subfigure}	
	\begin{subfigure}[b]{0.3\textwidth}
	\centering
    		\includegraphics[scale=1.3]{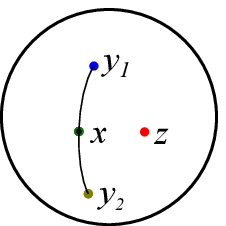}
		\caption{}
	\end{subfigure}
\caption{}
\label{fig:factoring}
\end{figure}

Рассмотрим, к примеру, карту $\M$, симметричную относительно поворота квадрата, представляющего тор, на угол в $180^\circ$ (fig. \ref{fig:factoring}). В этом случае мы вместо карты $\M$ можем рассматривать некоторую фактор-карту $\mathfrak{M}$, нарисованную  на сфере. Это была бы в точности карта на сфере с количеством вершин, ребер и граней, в два раза меньшим числа вершин, ребер и граней исходной карты $\M$, если бы у нас не было точек ветвления. Наличие точек ветвления может увеличить нам количество вершин, ребер и граней фактор-карты $\mathfrak{M}$ на сфере. Так, если какая-то вершина на фактор-карте $\mathfrak{M}$ не попадает ни в одну из особых точек, то она распадется на $L=2$ вершины исходной карты $\M$ на торе. Если же она, например, попадет в центр квадрата, то есть в особую точку, то тогда ей на торе будет соответствовать единственная вершина на исходной карте.  

В более общем случае предположим, что одна из вершин фактор-карты $\mathfrak{M}$ совпадает с точкой ветвления индекса $m_i$ орбифолда $O$. Тогда такой вершине будет отвечать $L/m_i$ вершин карты $\M$ на исходном многообразии $X$. При этом степень $d$ этой вершины $\mathfrak{x}$ увеличится в $m_i$ раз и окажется равной $m_id$. Действительно, все точки, принадлежащие окрестности точки ветвления, включая и те, которые принадлежат исходящим из вершины $\mathfrak{x}$ ребрам, особыми не являются, а следовательно, должны размножиться в $L$ раз. Для того, чтобы это обеспечить, мы и должны любой луч, исходящий из $\mathfrak{x}$, а также любой сектор, ограниченный соседними лучами, размножить в $m_i$ раз. 

Предположим теперь, что точка ветвления орбифолда $O$ попала в центр грани $f$ фактор-карты $\mathfrak{M}$ на орбифолде. Эта точка размножится в $L/m_i$ раз, где $m_i$ --- индекс точки ветвления. Остальные точки грани особыми не являются. Поэтому, как и в случае, когда точка ветвления попадает в вершину, ``степень'' грани, то есть количество инцидентных ей ребер и вершин, должна при поднятии на многообразие $X$ увеличиться в $m_i$ раз.  

Осталось разобраться с ребрами и понять, откуда у фактор-карты $\mathfrak{M}$ появляются висячие полуребра. Для этого рассмотрим вначале простейшую четырехвалентную карту $\M$ на торе $X$, состоящую из единственной вершины, двух ребер и одной грани (рис. \ref{fig:factoring}(a)). На квадрате, представляющем тор, такую карту мы можем нарисовать, поместив единственную вершину в центр квадрата, а два ребра нарисовать так, чтобы они пересекали особые точки $b_1$ и $b_2$ квадрата (рис. \ref{fig:factoring}(b)). Эта карта переходит в себя при поворотах на угол в $180^\circ$. Как следствие, на орбифолде $O$, представляющем собой сферу с четырьмя точками ветвления индексов $m_i=2$, карте $\M$ будет отвечать некоторая фактор-карта $\mathfrak{M}$, у которой имеется единственная вершина $\mathfrak{x}$ степени $2$ и единственная грань $f$. Однако из вершины $\mathfrak{x}$ фактор-карты будут исходить уже не ребра, а висячие полуребра $d_1$ и $d_2$, начинающиеся в вершине $\mathfrak{x}$ и заканчивающиеся в двух точках ветвления орбифолда $O$ индекса $2$ (рис. \ref{fig:factoring}(c)). Каждое из этих полуребер при поднятии на накрывающее $O$ многообразие $X$ превратится в нормальное ребро карты $\M$. 

Аналогичная ситуация имеет место и в общем случае. Именно, любая фактор-карта $\mathfrak{M}$ может состоять из ребер, свободных от точек ветвления, либо из полуребер, заканчивающихся в точках ветвления орбифолда $O$ индексов $m_i=2$. Каждое ребро $e$ фактор-карты $\mathfrak{M}$ размножится в $L$ ребер карты $\M$ на накрывающем $O$ многообразии $X$, а каждое висячее полуребро $d$ фактор-карты $\mathfrak{M}$ размножится в $L/2$ ребер карты $\M$. 

Заметим, что несколько различных гомеоморфизмов могут отвечать одному и тому же орбифолду. Примером служат два вращения квадрата, отождествляемого с тором, а именно, вращения на $90^\circ$ и на $270^\circ$. В работе \cite{Mednykh_Nedela} такие гомеоморфизмы периода $L$ отождествляются с эпиморфизмами группы автоморфизмов поверхности рода $g$ в циклическую группу $\Z_L$. Там же доказано, что общее число этих эпиморфизмов расcчитывается по формуле  
\begin{equation}
\label{eq:signat_number}
\Epi_0(\pi_1(O),\Z_L)=m^{2\mathfrak{g}}\phi_{2\mathfrak{g}}(L/m)\,E(m_1,\ldots,m_r),
\end{equation}
где $m=\lcm(m_1,\ldots,m_r)$ --- наименьшее общее кратное индексов ветвления,
$$
E(m_1,\ldots,m_r)=\dfrac{1}{m}\sum\limits_{k=1}^m\prod\limits_{i=1}^r\Phi(k,m_i),\qquad \text{где}\qquad
\Phi(k,m_i)=\dfrac{\phi(m_i)}{\phi(n_i)}\,\mu(n_i),\qquad n_i=\dfrac{m_i}{\gcd(k,m_i)},
$$
$\phi(n)$ --- функция Эйлера, $\mu(n)$ --- функция Мебиуса. При этом нам следует также учитывать дополнительные ограничения на индексы $m_1,\ldots,m_r$, основным из которых является формула Римана-Гурвица
\begin{equation}
\label{eq:Riemann-Hurvitz}
2-2g=L\left(2-2\mathfrak{g}-\sum\limits_{i=1}^r\left(1-\dfrac{1}{m_i}\right)\right).
\end{equation}
а также равенства вида
\begin{equation}
\label{eq:H2}
m=\lcm(m_1,\ldots,m_r)\mid L; m=L\,\, \text{где} \,\, g=0 ;
\end{equation}
\begin{equation}
\label{eq:H3}
\lcm(m_1,\ldots,m_{i-1},m_{i+1},\ldots,m_r)=m\quad\forall\,\,i=1,\ldots,r.
\end{equation}
Помимо данного равенства, нам в данной работе понадобятся следующие ограничения на числа $m_i$ \cite{Mednykh_Nedela}:
\begin{itemize}
\item если $m$ четно, количество $m_i$, делящихся на наибольшую степень двойки, которая делит $m$, четно;
\item если $G \geq 2$, то $r \neq 1$ и $r \geq 3$ for $g=0$, если $G=1$, то $r \in \{0,3,4\}$; наконец, если $G=0$, то $r=2$;
\item если $g >1$, то $L \leq 4g +2$.
\end{itemize}

Подводя итоги вышесказанному, мы можем из леммы Берсайда получить следующую формулу для подсчета количества непомеченных карт на многообразии рода $g$, обладающих теми или иными специальными свойствами \cite[Th. 3.1]{Mednykh_Nedela}:
\begin{equation}
\label{eq:mednyh_main}
\tilde{\epsilon}(g) = \frac{1}{2n}\sum_{L|2n} \sum_{O \in \Orb(S_g, \Z_L)} \Epi_0(\pi_1(O),\Z_L) \epsilon_O(2n/L)
\end{equation}
Здесь $\epsilon_O(2n/L)$ есть количество корневых фактор-карт с $2n/L$ полуребрами на орбифолде $O$, отвечающих картам на поверхности рода $g$ с заданными свойствами, в нашем случае – $d$-валентным однофейсовым картам.

Чтобы воспользоваться этой формулой, нам остается решить две подзадачи. Во-первых, учитывая ограничения (\ref{eq:Riemann-Hurvitz}), (\ref{eq:H2}) и (\ref{eq:H3}), описать все реализующиеся при этих ограничениях орбифолды. Во-вторых, по данному орбифолду понять, как устроены на нем карты, отвечающие $d$-валентным однофейсовым картам на поверхности рода $g$, и перечислить их.  

\section{Перечисление регулярных карт с одной гранью}

Рассмотрим, прежде всего, некоторые особенности перечисления произвольных однофейсовых регулярных карт. Если исходная карта $\M$ на поверхности $X$ была однофейсовой, то таковой будет и фактор-карта $\mathfrak{M}$ на орбифолде $O$. Но для того, чтобы единственная грань $\mathfrak{M}$ на орбифолде $O$ затем не размножилась, в нее обязательно нужно поместить одну из точек ветвления с индексом $m=L$. Иными словами, мы должны рассматривать только такие орбифолды $O$, в которых хотя бы одна точка ветвления обязательно присутствует, причем индекс этой точки ветвления обязан быть равным $L$. Кроме того, из условия (\ref{eq:H2}) следует, что 
$$
m=\lcm(m_1,\ldots,m_{r-1},L)\mid L\qquad \Longrightarrow\qquad m=L.
$$

Теперь предположим, что все вершины однофейсовой карты $\M$ имеют одну и ту же степень $d$. Заметим, что все остальные точки ветвления орбифолда $O$ должны попадать либо в вершины, либо совпадать с одним из двух концов полуребер фактор-карты $\mathfrak{M}$. При этом, как уже было отмечено выше, полуребра $\mathfrak{M}$ могут оканчиваться только в точках ветвления индекса $m_i=2$. Точка же ветвления индекса $m_i\neq 1$ может совпадать с вершиной $\mathfrak{x}$ фактор-карты $\mathfrak{M}$ только лишь в том случае, если $m_i$ является делителем степени $d$ вершин карты $\M$. При этом степень $\mathfrak{x}$ такой вершины равна $d/m_i$. 

Мы установили выше, что $m=L$. Выбирая в условии (\ref{eq:H3}) в качестве $i$-й точки ветвления точку, попавшую в грань фактор-карты $\mathfrak{M}$, мы получаем, что
\begin{equation}
\label{eq:H3_1}
m=L=\lcm(m_1,\ldots,m_{r-1}).
\end{equation}
В случае, когда $d$ есть четное число, все $m_i$, $i=1,\ldots,r-1$, делят $d$, а когда $d$ есть нечетное число, все такие $m_i$ делят $2d$. Условие (\ref{eq:H3_1}) тогда означает, что и период автоморфизма $L$ также должен быть либо делителем $d$ в случае $d$ четного, либо делителем $2d$ в случае $d$ нечетного.

Рассмотрим теперь случай простого $d>2$. В этом случае максимальное значение периода $L=2d$, а все возможные значения этого параметра равны $1,2,d,2d$. Значение $L=1$ соответствует всем корневым $d$-регулярным картам, построенным на поверхности рода $g$. Используя результаты работы \cite{Walsh_Lehman}, можно показать, что количество таких карт рассчитывается по формуле
\begin{equation}
\epsilon^{(d)}_\mathfrak{g}(k) = \frac{2(k+2g-1)!}{2^{2g}\cdot k!} \sum_{
\substack{i_1+\dots+i_k = g \\
i_1,\dots, i_k \geq 0}
}\prod_{j=1}^{k} \frac{\binom{d-1}{2i_j}}{2i_j+1}.
\label{one_face_labelled}
\end{equation}
Нам нужно описать слагаемые в формуле (\ref{eq:mednyh_main}), отвечающие оставшимся значениям $L$. Для этого рассмотрим вначале строение соответствующих этим $L$ орбифолдов $O$.

Из соотношений (\ref{eq:H3_1}) следует, что в случае $L=2$ все точки ветвления имеют один и тот же индекс ветвления, равный двум. С учетом этого наблюдения уравнение Римана-Гурвица дает следующее выражение, связывающее количество $r$ точек ветвления с параметрами $g$ и $\mathfrak{g}$:  
$$
r=2g+2-4\mathfrak{g},
$$
При этом количество эпиморфизмов (\ref{eq:signat_number}) оказывается равным $4^\mathfrak{g}$.
 
В случае $L=d$ все точки ветвления имеют индексы $m_i=d$. Наименьшее общее кратное $m$ этих чисел также равно $d$, а уравнение Римана-Гурвица принимает вид
$$
2-2g=d(2-2\mathfrak{g}-\sum\limits_{i=1}^{r_d}(1-1/d))=2d-2d\mathfrak{g}-(d-1)r_d,
$$
где $r_d$ --- количество точек ветвления индекса $d$. Из этого уравнения следует, что $\mathfrak{g}$ меняется от нуля до максимального значения
$$
\mathfrak{g}^{\max}=\left\lfloor\dfrac{2g+d-1}{2d} \right\rfloor,
$$
соответствующего $r_d^{\min}=1$. Несложно показать, что в данном случае количество эпиморфизмов рассчитывается по формуле
$$
\Epi_0(\pi_1(O),\Z_d)=d^{2\mathfrak{g}-1}\left[(d-1)(-1)^{r_d}+(d-1)^{r_d}\right].
$$

Теперь перейдем к случаю $L=2d$. Для такого $L$ единственная точка ветвления с индексом $m_i=2d$ принадлежит грани фактор-карты на орбифолде, $r_2$ точек ветвления, имеющие индексы, равные двум, отвечают полуребрам фактор-карты $\mathfrak{M}$, а $r_d$ точек ветвления индексов $d$ соответствуют вершинам этой фактор-карты $\mathfrak{M}$ степени один. С учетом приведенных выше обозначений уравнение Римана-Гурвица принимает вид
$$
2-2g=2d\left(2-2\mathfrak{g}-\dfrac{r_2}{2}-\dfrac{(d-1)}{d}r_d-\dfrac{(2d-1)}{2d}\right)
\quad \Longrightarrow\quad 2g+2d-1=4\mathfrak{g}d+dr_2+2(d-1)r_d.
$$
Из уравнения видно, что $r_2$ обязано быть нечетным числом (и, как следствие, обязано быть больше нуля). Кроме того, $r_d$ также должно быть больше нуля --- если мы исключим точку ветвления с индексом $m_i=2d$, то из условия
$$
\lcm(m_1,\ldots,m_{i-1},m_{i+1},\ldots,m_r)=2d
$$
мы получаем, что в этом наборе обязан присутствовать хотя бы один индекс $m_i=d$. Количество же эпиморфизмов (\ref{eq:signat_number}) в данном случае равно
$$
\Epi_0(\pi_1(O),\Z_{2d})=2\cdot(2d)^{2\mathfrak{g}-1}\cdot (d-1)\left[(d-1)^{r_d}-(-1)^{r_d}\right].
$$

Перейдем теперь ко второму этапу --- описанию фактор-карт на соответствующих орбифолдах. Исходная карта $\M$ на поверхности $X$ рода $g$ имеет $(4g-2)/(d-2)$ вершин степени $d$ и $(2g-1)d/(d-1)$ ребер. Опишем характеристики соответствующих $\M$ фактор-карт $\mathfrak{M}$ на орбифолде рода $\mathfrak{g}$ для различных значений параметра $L$.

Рассмотрим вначале фактор-карту $\mathfrak{M}$, отвечающую $L=2$. На соответствующем этому значению $L$ орбифолде $O$ имеется $r$ точек ветвления индекса $2$. Одна из этих точек гарантированно попадает на единственную грань фактор-карты $\mathfrak{M}$. Остальные точки могут попадать только в полуребра. Вершины фактор-карты $\mathfrak{M}$ имеют степень, равную $d$. Удобно, как обычно, считать, что полуребра фактор-карты $\mathfrak{M}$ --- это обычные ребра некоторой обычной же карты $\tilde{\M}$ на поверхности рода $\mathfrak{g}$ с листьями на конце. Обозначим через $s$ количество листьев этой карты. Мы только что установили, что $s=r-1$. Количество $k$ $d$-валентных вершин на этой карте $\tilde{\M}$ в $L=2$ раза меньше количества вершин на исходной карте, то есть равняется 
$$
k=\dfrac{2g-1}{d-2}.
$$
Наконец, количество полуребер у карты  $\tilde{\M}$ равно количеству ребер исходной карты плюс количеству $s=r_2-1$ полуребер, ведущих в листья этой карты $\tilde{\M}$:
$$
2n=\dfrac{d(2g-1)}{d-2}+r-1=\dfrac{d(2g-1)}{d-2}+2g-4\mathfrak{g}+1.
$$

Перейдем теперь к описанию соответствующей $\M$ фактор-карты $\mathfrak{M}$ в случае $L=d$. При таком значении параметра $L$ у нас одна точка ветвления индекса $d$ обязательно попадает в грань, а оставшиеся $s=r_d-1$ точек ветвления индекса $d$ должны совпадать с одновалентными вершинами фактор-карты $\mathfrak{M}$. При этом количество $k$ $d$-валентных вершин у нас в данном случае оказывается равным
$$
k=\dfrac{V-s}{d}=\dfrac{4g-2}{d(d-2)}-\dfrac{r_d-1}{d},
$$
а количество $n$ ребер уменьшается по сравнению с количеством ребер исходной карты $\M$ в $d$ раз:
$$
n=\dfrac{2g-1}{d-2}.
$$

Наконец, в случае $L=2d$ у нас, опять-таки, одна точка ветвления индекса $2d$ попадает в грань, $r_d$ точек ветвления индекса $d$ попадают в вершины фактор-карты $\mathfrak{M}$ степени $1$, а $r_2$ точек ветвления индекса $2$ попадают в полуребра этой фактор-карты, превращаясь в листья обычной карты $\tilde{\M}$. Таким образом, количество $s$ листьев карты $\tilde{\M}$ равно $r_2+r_d$. Для подсчета количества $k$ трехвалентных вершин нам из числа вершин исходной карты $\M$ нужно отнять удвоенное количество $2r_d$ вершин, в которые попали точки ветвления индекса $d$, а затем поделить результат на $L=2d$: 
$$
k=\dfrac{2g-1}{d(d-2)}-\dfrac{r_d}{d}.
$$
Количество же полуребер $2n$ у такой карты $\tilde{\M}$ равняется количеству полуребер исходной карты, деленному на $L=2d$, плюс количеству полуребер $r_2$, ведущих в листья, совпадающие с точками ветвления индекса $2$:
$$
2n=\dfrac{2g-1}{d-2}+r_2.
$$ 

Перейдем теперь к непосредственному подсчету некорневых $d$-валентных карт с одной гранью в случае простого $d>2$. Обозначим через  $\epsilon^{(1\div d)}_\mathfrak{g}(k)$ количество корневых $(1\div d)$-валентных карт $\tilde{\M}$ с одной гранью на поверхности рода $\mathfrak{g}$, имеющих $k$ вершин степени $d$. В случае $d=3$ у нас имеется аналитическое выражение \cite{Chapuy} для этих чисел: 
$$
\epsilon^{(1\div 3)}_\mathfrak{g}(k)=\dfrac{2\,n!}{12^\mathfrak{g}\,\mathfrak{g}!\,(k-\mathfrak{g})!\,s!}
=\dfrac{2(2k-2\mathfrak{g}+1)!}{12^\mathfrak{g}\mathfrak{g}!(k-\mathfrak{g})!(k+2-4\mathfrak{g})!}.
$$
В остальных случаях можно использовать формулу, полученную в работе Уолша и Лемана \cite{Walsh_Lehman}: 
\begin{equation}
\epsilon^{(1 \div d)}_\mathfrak{g}(k) = \frac{2\cdot (d\cdot k-k-2g+1)!}{2^{2g}\cdot k! \cdot (2-4g+(d-2)k)!} \sum_{
\substack{i_1+\dots+i_k = g \\
i_1,\dots, i_k \geq 0}
}\prod_{j=1}^{k} \frac{\binom{d-1}{2i_j}}{2i_j+1}.
\label{one_face_labelled}
\end{equation}

В случае $L=2$ у исходной $1\div d$-карты $\tilde{\M}$ имеется $2n=d(2g-1)/(d-2)+s$ способов выбрать одно из полуребер корневым. Однако $s$ из них отвечают полуребрам фактор-карты $\mathfrak{M}$ на орбифолде, и потому корневыми быть не могут. Используя принцип двойного подсчета, мы должны для определения количества $1\div d$-карт $\mathfrak{M}$ на орбифолде домножить $\epsilon^{(1\div 3)}_\mathfrak{g}(k)$ на число возможных полуребер фактор-карты $\mathfrak{M}$ и поделить на количество $2n$ полуребер карты $\tilde{\M}$: 
$$
\dfrac{d(2g-1)}{2n(d-2)}\cdot\epsilon^{(1\div d)}_\mathfrak{g}(k).
$$
Домножая полученное выражение на количество $4^{\mathfrak{g}}$ эпиморфизмов и суммируя по $\mathfrak{g}$ от $0$ до $\lfloor g/2\rfloor$, мы для случая $L=2$ получим выражение вида
$$
f_2(g):=\dfrac{d(2g-1)}{d-2}\sum\limits_{\mathfrak{g}=0}^{\lfloor g/2\rfloor}
\dfrac{4^\mathfrak{g}}{2n(\mathfrak{g})}\cdot \epsilon^{(1\div d)}_\mathfrak{g}(k(\mathfrak{g})).
$$
В частности, для $3$-валентных карт имеем следующее явное аналитическое выражение:
$$
f_2(g)=3(2g-1)\sum\limits_{\mathfrak{g}=0}^{\lfloor g/2\rfloor}
\dfrac{(n-1)!}{3^\mathfrak{g}\,\mathfrak{g}!\,(k-\mathfrak{g})!\,s!}=
(2g-1)\sum\limits_{\mathfrak{g}=0}^{\lfloor g/2\rfloor}
\dfrac{(4g-2-2\mathfrak{g})!}{3^{\mathfrak{g}-1}\,\mathfrak{g}!\,(2g-1-\mathfrak{g})!\,(2g-4\mathfrak{g}+1)!}.
$$

В случае $L=d$ нам нужно просто домножить количество $\epsilon^{(1\div d)}_\mathfrak{g}(k)$ корневых $(1\div d)$-карт на соответствующее данному $L$ количество эпиморфизмов и просуммировать результат по $\mathfrak{g}$ от нуля до $\mathfrak{g}^{\max}=\lfloor(2g+d-1)/(2d) \rfloor$:
$$
f_d(g):=\sum\limits_{\mathfrak{g}=0}^{\lfloor(2g+d-1)/(2d) \rfloor}d^{2\mathfrak{g}-1}\left[(d-1)(-1)^{r_d(\mathfrak{g})}+(d-1)^{r_d(\mathfrak{g})}\right]\cdot 
\epsilon^{(1\div d)}_\mathfrak{g}(k(\mathfrak{g})).
$$
Для трехвалентных карт отсюда получаем, что
$$
f_3(g)=\sum\limits_{\mathfrak{g}=0}^{\lfloor(g+1)/3 \rfloor}\dfrac{2\,n!\,3^{2\mathfrak{g}-1}\left[2^{r_3}+2(-1)^{r_3}\right]}{12^\mathfrak{g}\,\mathfrak{g}!\,(k-\mathfrak{g})!\,s!}=
\dfrac{(2g-1)!}{(g-1)!}\sum\limits_{\mathfrak{g}=0}^{\lfloor(g+1)/3 \rfloor}\left(\dfrac{3}{4}\right)^{\mathfrak{g}-1}
\dfrac{2^{g+1-3\mathfrak{g}}+(-1)^{g+2-3\mathfrak{g}}}{\mathfrak{g}!\,(g+1-3\mathfrak{g})!}.
$$

Наконец, в случае $L=2d$ нам, как и для значения $L=2$, нужно в процессе перехода от $(1\div d)$-карты $\tilde{\M}$ к фактор-карте $\mathfrak{M}$ на орбифолде пересчитать количество возможных корневых полуребер, домножая $\epsilon^{(1\div d)}_\mathfrak{g}(k)$ на число $(2g-1)/(d-2)$ полуребер фактор-карты $\mathfrak{M}$ и деля на количество $2n$ полуребер $(1\div d)$-карты $\tilde{\M}$. Помимо этого, нам в полученной фактор-карте $\mathfrak{M}$ нужно $\BCf{s}{r_2}$ способами выбрать из $s=r_d+r_2$ листьев те, которые попадают в $r_2$ точек ветвления орбифолда индекса $2$. Наконец, полученное выражение нам нужно домножить на количество эпиморфизмов, а затем просуммировать его по параметрам $\mathfrak{g}$ и $k$:
$$
f_{2d}(g):=\dfrac{2(2g-1)}{d-2}\sum\limits_{\mathfrak{g},k} \dfrac{(2d)^{2\mathfrak{g}-1}\cdot (d-1)\left[(d-1)^{r_d}-(-1)^{r_d}\right]}{2n}
\cdot \BCf{r_2+r_d}{r_2}\cdot \epsilon^{(1\div d)}_\mathfrak{g}(k).
$$
Для случая трехвалентных карт получаем отсюда, что
$$
f_6(g)=2(2g-1)\sum\limits_{\mathfrak{g},k}\dfrac{3^{\mathfrak{g}-1}\left[2^{r_3}-(-1)^{r_3}\right]\,(n-1)!}
{\mathfrak{g}!\,(k-\mathfrak{g})!\,r_2!\,r_3!}=
$$ 
$$
=2(2g-1)\sum\limits_{k=k_{\min}}^{k_{\max}}\sum\limits_{\mathfrak{g}=0}^{\mathfrak{g}_{\max}(k)}
\dfrac{3^{\mathfrak{g}-1}\left[2^{2g-1-3k}-(-1)^{2g-1-3k}\right]\,(2k-2\mathfrak{g})!}
{\mathfrak{g}!\,(k-\mathfrak{g})!\,(4k+3-2g-4\mathfrak{g})!\,(2g-1-3k)!}.
$$
Граничные значения для $k$ и $\mathfrak{g}$ получаются из формул
$$
r_3=2g-1-3k,\qquad \qquad r_2=4k+3-2g-4\mathfrak{g},
$$
полученных из формулы Римана-Гурвица и выражения для $k$. Так, максимальное значение $k$ достигается при минимальном значении $r_3=1$ и равно
$$
k_{\max}=\left\lfloor \dfrac{2g-2}{3}\right\rfloor.
$$
Минимальное же значение $k$ достигается при $\mathfrak{g}=0$ и $r_2=r_2^{\min}$:
$$
k_{\min}=\left\lfloor \dfrac{2g-3+r_2^{\min}}{4}\right\rfloor.
$$
При четном $g$ значение $r_2^{\min}=3$, а при нечетном оно равно единице. В обоих случаях
$$
k_{\min}=\left\lfloor \dfrac{g}{2}\right\rfloor.
$$
Максимальное же значение $\mathfrak{g}$ при фиксированных $g$ и $k$ достигается в случае, когда $r_2=r_2^{\min}$ и равно
$$
\mathfrak{g}_{\max}=\left\lfloor \dfrac{4k-2g+3-r_2^{\min}}{4}\right\rfloor=k-\left\lfloor \dfrac{g}{2}\right\rfloor.
$$ 
Таким образом,
$$
f_6(g)=2(2g-1)\sum\limits_{k=\lfloor g/2\rfloor}^{\lfloor (2g-2)/3\rfloor}\sum\limits_{\mathfrak{g}=0}^{k-\lfloor g/2\rfloor}
\dfrac{3^{\mathfrak{g}-1}\left[2^{2g-1-3k}-(-1)^{2g-1-3k}\right]\,(2k-2\mathfrak{g})!}
{\mathfrak{g}!\,(k-\mathfrak{g})!\,(4k+3-2g-4\mathfrak{g})!\,(2g-1-3k)!}.
$$

Подставляя теперь выражения для $\epsilon^{(3)}(g)$, $f_2(g)$, $f_3(g)$ и $f_6(g)$ в (\ref{eq:mednyh_main}), 
получаем следующую окончательную формулу для подсчета количества всех трехвалентных непомеченных карт с одной гранью на многообразии рода $g$:
$$
\tilde{\epsilon}^{(3)}(g)=\dfrac{(6g-4)!}{12^g\,g!\,(3g-2)!}+
\sum\limits_{\mathfrak{g}=0}^{\lfloor g/2\rfloor}
\dfrac{(4g-2-2\mathfrak{g})!}{2\cdot 3^{\mathfrak{g}}\,\mathfrak{g}!\,(2g-1-\mathfrak{g})!\,(2g-4\mathfrak{g}+1)!}+
$$
\begin{equation}
\label{eq:num_nonrooted_3}
+\dfrac{(2g-2)!}{6\cdot (g-1)!}\sum\limits_{\mathfrak{g}=0}^{\lfloor(g+1)/3 \rfloor}\left(\dfrac{3}{4}\right)^{\mathfrak{g}-1}
\dfrac{2^{g+1-3\mathfrak{g}}+(-1)^{g-\mathfrak{g}}}{\mathfrak{g}!\,(g+1-3\mathfrak{g})!}+
\end{equation}
$$
+\sum\limits_{k=\lfloor g/2\rfloor}^{\lfloor (2g-2)/3\rfloor}\sum\limits_{\mathfrak{g}=0}^{k-\lfloor g/2\rfloor}
\dfrac{3^{\mathfrak{g}-2}\left[2^{2g-1-3k}+(-1)^{k}\right]\,(2k-2\mathfrak{g})!}
{\mathfrak{g}!\,(k-\mathfrak{g})!\,(4k+3-2g-4\mathfrak{g})!\,(2g-1-3k)!}.
$$

Описанный выше подход можно использовать для получения аналогичных формул и для случая, когда $d$ не является простым числом. Так, в случае $d=4$ с использованием похожих рассуждений получается следующая формула для перечисления $4$-регулярных карт с одной гранью:
$$
\tilde{\epsilon}^{(4)}(g)=\dfrac{(4g-3)!}{4^g\,g!\,(g-1)!}+\dfrac{3(4g-3)!}{2\,(2g+1)!\,(2g-2)!}+
$$
\begin{equation}
\label{eq:num_nonrooted_4}
+\sum\limits_{\mathfrak{g}=1}^{\lfloor g/2\rfloor}\sum\limits_{k=2\mathfrak{g}-1}^{g-1}\dfrac{(2g-2\mathfrak{g}+k-1)!}{2\,(2k-4\mathfrak{g}+2)!\,\mathfrak{g}!\,(k-\mathfrak{g})!\,(2g-1-2k)!}+
\end{equation}
$$
+\sum\limits_{\mathfrak{g}=0}^{\lfloor g/4\rfloor}\sum\limits_{\substack{r_4=2\\ 2 \mid r_4}}^{\lfloor2(g+3-4\mathfrak{g})/3\rfloor}
\sum\limits_{k=2\mathfrak{g}-1+r_4/2}^{\lfloor g/2+r_4/4\rfloor}
\dfrac{2^{2\mathfrak{g}-3+r_4}\,(k-2\mathfrak{g}+g-r_4/2)!}{\mathfrak{g}!\,(k-\mathfrak{g})!\,(g-r_4/2-2k)!\,(2k+3-4\mathfrak{g}-r_4)!\,(r_4-1)!}.
$$

Соответствующие числа приведены в таблицах (\ref{table_first_one_face}) --- (\ref{table_last_one_face}).

\section{Перечисление регулярных карт с несколькими гранями}

Для перечисления корневых регулярных карт с несколькими гранями мы будем пользоваться походом декомпозиции по корневом ребру (root-edge decomposition), предложенным Таттом \cite{Tutte_Census}. Его идея состояла в том, чтобы классифицировать карты в зависимости от типа структуры, получающейся в результате удаления корневого ребра, а затем подсчитать количество карт в каждом классе. Перечисление $d$-регулярных карт с помощью этого подхода имеет определенные особенности. Прежде всего, заметим, что при удалении корневого ребра могут изменится степени двух вершин, так что из регулярной карты мы получим некоторую другую карту, имеющую особые вершины. Мы можем уменьшить связанные с этим трудности, если вместо удаления ребра будем применять стягивание. Тем не менее, и при стягивании ребра степень одной из вершин будет изменяться.

\begin{figure}[ht]
\centering
\includegraphics[scale=0.9]{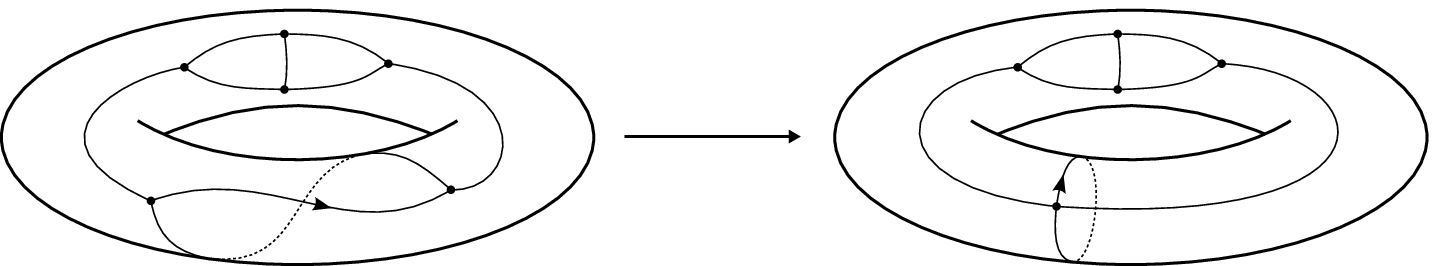}	
\caption{}
\label{fig:case1}
\end{figure}

Рассмотрим для примера $d$-регулярную карту на торе, изображенную на рис. \ref{fig:case1}. Ее корневое ребро соединяет две различные вершины, при этом его стягивание приводит к появлению новой карты на поверхности того же рода (рис \ref{fig:case1}), корень у которой можно выбрать неким каноническим способом. Как видно, степень новой корневой вершины отлична от $d$. Возможны и более сложные случаи.

\begin{figure}[ht]
\centering
\includegraphics[scale=0.9]{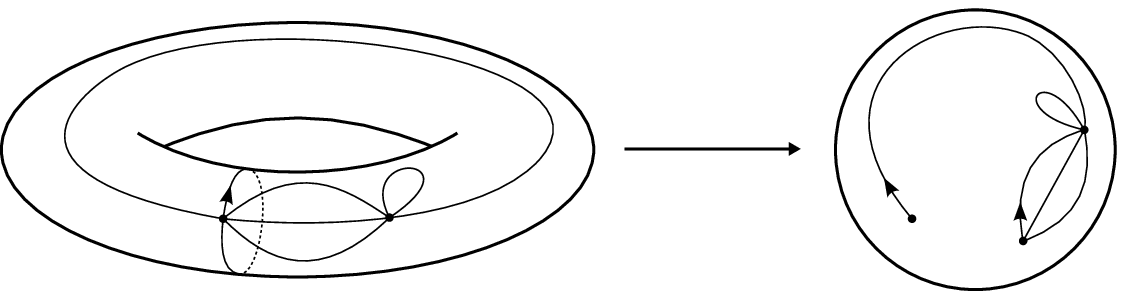}	
\caption{}
\label{fig:case2}
\end{figure}

Действительно, рассмотрим карту, изображенную на \ref{fig:case2}. Ее корневое ребро представляет собой петлю, а значит, стянуть его как раньше уже нельзя. Поэтому мы в данном случае будем понимать под стягиванием некоторую более сложную операцию. Она будет заключаться в разрезании поверхности вдоль корневого ребра и стягивании получившиеся краев в точки. Можно думать, что эти две получившиеся точки есть результат разделения корневой вершины (той, которой инцидентно корневое ребро) на две новых вершины, при каждой из которых каноническим образом выбирается корневое ребро. Если род поверхности достаточно высок, то, применяя подобную операцию несколько раз, мы получим карту с тремя и более корневыми ребрами, исходящими из различных вершин. Таким образом, для решения задачи о перечислении $d$-регулярных карт нам придется перечислить несколько более общий класс карт, а именно, карты, имеющие несколько корневых вершин заданных степеней помимо вершин степени $d$. Каждая такая корневая вершина будет иметь выделенное индицентное ей ребро. Карты описанного типа мы будем называть {\em почти регулярными}.

Пусть $Q(d;g;n;d_1,d_2,\ldots,d_k)$ есть число почти регулярных карт рода $g$, имеющих $k$ корневых вершин $v_1,v_2,\ldots,v_k$ степеней $d_1,d_2,\ldots,d_k$ соответственно и $n$ других некорневых вершин степени $d$. Мы выведем для этих чисел рекуррентное соотношение, классифицировав карты по типу структур, получающихся в результате стягивания первого корневого ребра. Все карты естественным образом разделяются в этом смысле на две класса. К первому классу относятся карты, в которых первое корневое ребро соединяет различные вершины. Второй класс состоит из карт, имеющих в качестве первого корневого ребра петлю.

Первый класс, в свою очередь, можно разделить на два непересекающихся подкласса. К первому подклассу отнесем карты, в которых корневое ребро соединяет корневую (помеченную) вершину с непомеченной. При стягивании такого ребра мы получим карту с числом непомеченных вершин, уменьшившимся на единицу, и степенью первой корневой вершины, увеличившейся на $d-2$ (рис. \ref{fig:case1}), так что число карт в данном подклассе равно $Q(d;g;n-1;d_1+d-2,d_2,\ldots,d_k)$.

\begin{figure}[ht]
\centering
\includegraphics[scale=0.9]{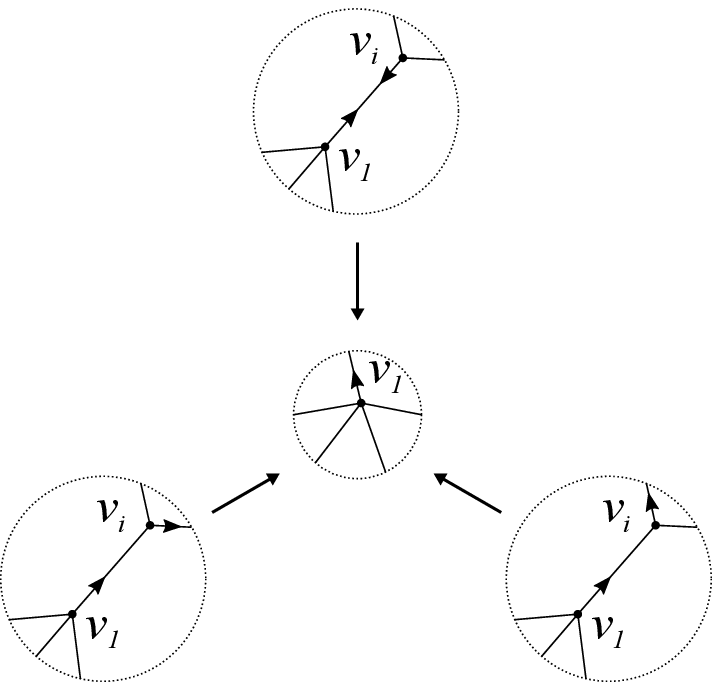}	
\caption{}
\label{fig:case1_2}
\end{figure}

Ко второму подклассу отнесем карты, в которых первое корневое ребро соединяет две различные выделенные вершины. При стягивании такого ребра две корневые вершины $v_1$ и $v_i$ стянутся в одну новую вершину, степень которой будет равна сумме степеней $v_1$ и $v_i$, уменьшенной на 2 (рис. \ref{fig:case1_2}). При новой вершине корневое ребро можно выбрать каноническим образом, определяемым корневым ребром при $v_1$. Для примера, можно выбрать ребро, идущее следующим в порядке обхода грани, находящейся слева от первого корневого ребра, против часовой стрелки. Заметим, что при этом мы потеряем информацию о том, какое из полуребер при $v_i$ было корневым, а значит, $d_i$ различных карт перейдут в одну и ту же карту под действием такой операции (рис \ref{fig:case1_2}). Как следствие, число карт в рассматриваемом подклассе равно такой сумме:
$$
\sum_{i=2}^k {d_i \cdot Q(d;g;n;d_1+d_i-2,d_2,\ldots,d_{i-1}, d_{i+1}, \ldots,d_k)}.
$$

\begin{figure}[ht]
\centering
\includegraphics[scale=0.9]{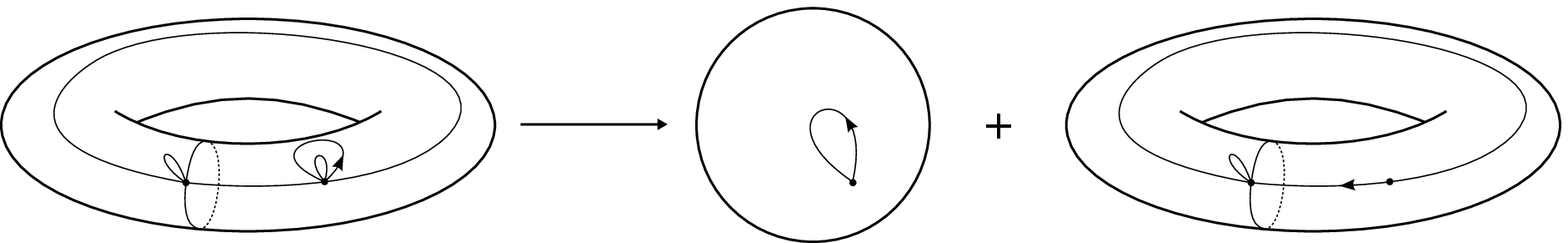}	
\caption{}
\label{fig:case3}
\end{figure}

Перейдем теперь к перечислению второго класса карт, а именно, карт, в которых первое корневое ребро является петлей. Разрезая поверхность вдоль этого ребра и стягивая получившиеся края, мы получим поверхность без края, которая либо связна и имеет род $g-1$ (рис \ref{fig:case2}), либо состоит из двух связных компонент, сумма родов которых равна роду исходной поверхности (рис. \ref{fig:case3}). В обоих случаях две точки, полученные стягиванием границ, есть результат разделения первой корневой вершины на две новых вершины, при каждой из которых можно канонически выбрать корневое ребро. Сумма степеней этих вершин на 2 меньше, чем степень исходной вершины.

В первом случае, случае уменьшения рода, мы немедленно получаем следующую перечислительную формулу:
$$
\sum_{d_1'=0}^{d_1-2}Q(d;g-1;n;d_1',d_2,\ldots,d_k, d_1-2-d_1').
$$

Во втором случае мы получаем две отдельные карты, и нам необходимо перебрать все возможные способы распределить корневые вершины между ними, а также учесть все возможные комбинации родов. Это приводит к следующей перечислительной формуле для данного случая:
$$\sum_{g'=0}^{g} \sum_{d_1'=0}^{d_1-2}  \sum_{i=0}^{n} \sum_{I \subseteq \{2,\dots,k\}}Q(d;g';i;d_1',I) \cdot Q(d;g-g';n-i;d_1-2-d_1,\{2,\dots,k\} \backslash I).
$$
Под $Q(d,g,n,d_1,I)$ мы здесь понимаем $Q(d,g,n,d_1,d_{i_1},\dots, d_{i_m})$, где $i_1,\dots,i_m$ есть элементы $I$, выписанные в порядке неувеличения.

Подводя итог вышесказанному, мы получаем следующую формулу для перечисления почти регулярных карт:
$$
Q(d;g;n;d_1,d_2,\ldots,d_k) = Q(d;g;n-1;d_1+d-2,d_2,\ldots,d_k) +
$$
$$
+ \sum_{i=2}^k {d_i \cdot Q(d;g;n;d_1+d_i-2,d_2,\ldots,d_{i-1}, d_{i+1}, \ldots,d_k)} +
$$
$$
+ \sum_{d_1'=0}^{d_1-2}Q(d;g-1;n;d_1',d_2,\ldots,d_k, d_1-2-d_1')+
$$
$$
+ \sum_{\substack{
	g'+g''= g \\
	d'_1+d''_1=d_1-2\\
	n'+n''=n\\
	I \subseteq \{2,\dots,k\} \\
	g',g'',d',d'',n',n''\geq 0	
}} Q(d;g';i';d_1',I) \cdot Q(d;g'';i'';d_1'', \{2,\dots,k\} \backslash I).
$$
Начальные и граничные условия такие:
$$
\qquad Q(d;g;n;d_1,\dots,d_i,\dots,d_k) = 0 \quad \text {если} \quad g < 0, \, n < 0 \quad \text {или} \quad d_i \leq 0, \quad \text {кроме} \quad Q(d;0;0;0) = 1.
$$

Для перечисления $d$-регулярных карт рода $g$ с $n$ вершинами можно использовать $Q(d;g;n-1;d)$.

Перейдем далее к перечислению некорневых карт с несколькими гранями. Пусть $Q_O(d;h)$ есть число фактор-карт с $h$ полуребрами на орбифолде $O$, соответствующих $d$-регулярным картам на исходной поверхности рода $g$. Формула (\ref{eq:mednyh_main}) может быть переписана для этого случая как
$$
\tilde{Q}(d;g;n) = \frac{1}{d \cdot n}\sum_{L|(d \cdot n)} \sum_{O \in \Orb(S_g,Z_L)}\Epi_0(\pi_1(O), \Z_L) \cdot Q_O(d;d \cdot n/L),
$$
где $\tilde{Q}(d;g;n)$ есть число некорневых  $d$-регулярных карт рода $g$, имеющих $n$ вершин. Наша конечная цель — выразить $Q_O(d;h)$ через числа $Q(d;g;n;d_1,d_2,\ldots,d_k)$ почти регулярных карт.

Пусть сигнатура орбифолда $O$ равна $ (\mathfrak{g},L,[m_1,\ldots,m_r])$ и пусть $r_i$ есть число точек ветвления индекса $i$ в этой сигнатуре. Для тривиальных орбифолдов $Q_O(d;h)$ равно числу $d$-регулярных карт с $h/d$ вершинами:
$$
Q_O(d;h) = Q(d;\mathfrak{g}; h/d-1; d).
$$

Уиверждается, что для других орбифолдов выполнено такое равенство для $Q_O(d;h)$:

$$
Q_O(d;h) = \sum_{\substack{s=0\\2|(h-s)}}^{r_2}Q(d;\mathfrak{g};h/d - 1;d,\underbrace{1,\dots,1}_{s \, \text{times}}) \cdot \frac{1}{s!} \cdot \binom{2-2\mathfrak{g}+(h-s)/2-h/d}{r_2-s,r_3,r_4,\ldots} + 
$$

\begin{equation}
+\sum_{\substack{
	0\leq v_i \leq r_i\\
	i \geq 2\\
	\exists i \colon v_i > 0\\
	i \nmid d \Rightarrow v_i = 0	
}} \sum_{\substack{s=0\\2|(h+s)}}^{r_2-v_2}Q \bigg( d;\mathfrak{g};\frac{h}{d} - \sum_{i=2}^d \frac{v_i}{i};
\underbrace{\frac{d}{2},\dots,\frac{d}{2}}_{v_2 \,\text{times}},
\underbrace{\frac{d}{3},\dots,\frac{d}{3}}_{v_3 \,\text{times}},
\dots,
\underbrace{\frac{d}{d},\dots,\frac{d}{d}}_{s+v_d \,\text{times}}
\bigg)\cdot 
\label{unlabelled_multiface}
\end{equation}
$$
\cdot \frac{h+s}{(v_d+s)! \cdot \prod_{i \neq d}(d/i)^{v_i} \cdot v_i!}\cdot
\binom{s+v_d}{s} \cdot \frac{h}{h+s}\cdot \binom{2-2\mathfrak{g}+(h-s)/2-h/d-\sum_{i=2}^d \frac{v_i(i-1)}{i}}{r_2-v_2-s,r_3-v_3,r_4-v_4,\ldots}.
$$

Это выражение для $Q_O(d;h)$ содержит два слагаемых, соответствующих фактор-картам без вершин, совпадающих с точками ветвления, и с такими вершинами. Начнем с первого слагаемого, которое описать легче. Оно состоит из суммы по $s$, числу точек ветвления, попадающих на висячие концы полуребер. Общее число $h$ полуребер равно  $2e+s$, где $e$ есть число ребер, так что $h-s$ обязательно будет четным: $2|(h-s)$. Рассматриваемые карты имеют $h/d$ вершин степени $d$, включая корень. Также они имеют $s$ висячих концов полуребер, на каждый из которых мы поместим по новой вершине, чтобы превратить фактор-карту в обыкновенную карту. Число таких карт дается выражением $Q(d;\mathfrak{g};h/d - 1;d,\underbrace{1,\dots,1}_{s \, \text{times}})$, которое дополнительно умножается на $1/s!$ чтобы не учитывать упорядоченность листьев. Мультиномиальный коэффициент в (\ref{unlabelled_multiface}) описывает число способов распределить оставшиеся точки ветвления орбифолда по гранями фактор-карты. Действительно, по формуле Эйлера число граней равно $2-2\mathfrak{g}+(h-s)/2-h/d$.

Во втором слагаемом в (\ref{unlabelled_multiface}) производятся достаточно похожие вычисления, однакого теперь мы должны учесть распределение точек ветвления по вершинам фактор-карты. Внешняя сумма берется по всем возможным способам распределить точки ветвления обрифолда по вершинам. Здесь $v_i$ есть число точек ветвления индекса $i$, совпадающих с вершинами. Как минимум одно из $v_i$ должно быть ненулевым, так как оставшиеся карты уже были перечислены в первом слагаемом. Также $v_i$ может быть ненулевым только если $i$ делит $d$. Внутренняя сумма берется по числу точек ветвления, совпадающих с висячими концами полуребер.

Для фиксированного распределения точек ветвления по вершинам, ребрам и граням, мы вначале вычисляем количество корневых карт рода $\mathfrak{g}$, имеющих соответствующее распределение степеней вершин. Далее, мы как и ранее заменяем висячие концы полуребер на листья, так что общее число листьев становится равно $s+v_d$.

Вспомним, что почти регулярные карты, перечисляемые $Q$, имеют упорядоченные корневые вершины, так что этот порядок необходимо убрать, разделив промежуточный результат на $(v_d+s)! \cdot \prod_{i \neq d}(d/i)^{v_i} \cdot v_i!$. Умножение на $h+s$ учитывает, что каждое полуребро может теперь быть корнем.

Биномиальный коэффициент $\binom{s+v_d}{s}$ подсчитывает различные способы выбрать $s$ из $s+v_d$ листьев для того, чтобы сделать из них висячие концы полуребер. Такая операция изменяет число способов выбрать в карте корень, так что мы должны воспользоваться двойным подсчетом и умножить промежуточный результат на $\frac{h}{h+s}$. Наконец, мультиномиальный коэффициент подсчитывает способы распределить точки ветвления, не совпадающие с вершинами и ребрами, по граням фактор-карты.

Формула (\ref{unlabelled_multiface}) была использована для перечисления $d$-регулярных карт рода $g$ для $g$ до 10 и для $d$ от $3$ до $6$. Соответствующие численные результаты приведены в таблицах (\ref{table_first}) --- (\ref{table_last}).

\newpage

\begin{table}[h!]
\begin{center}
\scriptsize
\begin{tabular}{c|ccc}
\midrule
$ g \backslash d$ &\phantom{00000000000000000}3\phantom{00000000000000000}&\phantom{00000000000000000}4\phantom{00000000000000000}&\phantom{0000000000000000}5\phantom{0000000000000000}\\ 
\midrule
1 & 1 & 1 & 0 \\
2 & 105 & 45 & 33 \\
3 & 50050 & 9450 & 0 \\
4 & 56581525 & 4729725 & 0 \\
5 & 117123756750 & 4341887550 & 1038647610 \\
6 & 386078943500250 & 6352181485650 & 0 \\
7 & 1857039718236202500 & 13566444744352500 & 0 \\
8 & 12277353837189093778125 & 39834473380605028125 & 3178849676735117385 \\
9 & 106815706684397824557193750 & 153946961458244898693750 & 0 \\
10 & 1183197582943074702620035168750 & 757572997336023146471943750 & 0 \\
11 & 16259070931137207808967206912537500 & 4625189759553876588251163487500 & 123443391148242936280071421860 \\
\midrule
\end{tabular}
\caption{Помеченные однофейсовые $d$-регулярные карты}
\label{table_first_one_face}
\end{center}
\end{table}

\begin{table}[h!]
\scriptsize
\centering
\begin{tabular}{c|ccc}
\midrule
$ g \backslash d$  &\phantom{00000000000000000}3\phantom{00000000000000000}&\phantom{00000000000000000}4\phantom{00000000000000000}&\phantom{0000000000000000}5\phantom{0000000000000000}\\ 
\midrule
1 & 1 & 1 & 0 \\
2 & 9 & 6 & 7 \\
3 & 1726 & 510 & 0 \\
4 & 1349005 & 169772 & 0 \\
5 & 2169056374 & 120644422 & 34629024 \\
6 & 5849686966988 & 144369379620 & 0 \\
7 & 23808202021448662 & 260893265836244 & 0 \\
8 & 136415042681045401661 & 663907896121296616 & 63576994019338897 \\
9 & 1047212810636411989605202 & 2263925904300525582790 & 0 \\
10 & 10378926166167927379808819918 & 9968065754464730977513732 & 0 \\
11 & 129040245485216017874985276329588 & 55061782851836038471634743076 & 1763477016403597971209426672 \\
\midrule
\end{tabular}
\caption{Непомеченные однофейсовые $d$-регулярные карты}
\label{table_last_one_face}
\end{table}

\begin{table}[h!]
\begin{center}
\scriptsize
\begin{tabular}{c|cccccc}
\midrule
$ v \backslash g$ &\phantom{0000000}0\phantom{0000000}&\phantom{0000000}1\phantom{0000000}&\phantom{0000000}2\phantom{0000000}&\phantom{0000000}3\phantom{0000000}&\phantom{0000000}4\phantom{0000000}&\phantom{0000000}5\phantom{0000000}\\ 
\midrule
2 & 4 & 1 & 0 & 0 & 0 & 0\\
4 & 32 & 28 & 0 & 0 & 0 & 0\\
6 & 336 & 664 & 105 & 0 & 0 & 0\\
8 & 4096 & 14912 & 8112 & 0 & 0 & 0\\
10 & 54912 & 326496 & 396792 & 50050 & 0 & 0\\
12 & 786432 & 7048192 & 15663360 & 6722816 & 0 & 0\\
14 & 11824384 & 150820608 & 544475232 & 518329776 & 56581525 & 0\\
16 & 184549376 & 3208396800 & 17388675072 & 30117189632 & 11100235520 & 0\\
18 & 2966845440 & 67968706048 & 522638463744 & 1465000951488 & 1191262520280 & 117123756750\\
20 & 48855252992 & 1435486650368 & 15007609257984 & 62975678300160 & 92809670660096 & 30625920998400\\
\midrule
\end{tabular}
\caption{Помеченные $3$-регулярные карты}
\label{table_first}
\end{center}
\end{table}

\begin{table}[h!]
\begin{center}
\scriptsize
\begin{tabular}{c|cccccc}
\midrule
$ v \backslash g$ &\phantom{0000000}0\phantom{0000000}&\phantom{0000000}1\phantom{0000000}&\phantom{0000000}2\phantom{0000000}&\phantom{0000000}3\phantom{0000000}&\phantom{0000000}4\phantom{0000000}&\phantom{0000000}5\phantom{0000000}\\ 
\midrule
2 & 2 & 1 & 0 & 0 & 0 & 0\\
4 & 6 & 5 & 0 & 0 & 0 & 0\\
6 & 26 & 46 & 9 & 0 & 0 & 0\\
8 & 191 & 669 & 368 & 0 & 0 & 0\\
10 & 1904 & 11096 & 13448 & 1726 & 0 & 0\\
12 & 22078 & 196888 & 436640 & 187580 & 0 & 0\\
14 & 282388 & 3596104 & 12974156 & 12350102 & 1349005 & 0\\
16 & 3848001 & 66867564 & 362330506 & 627525739 & 231290187 & 0\\
18 & 54953996 & 1258801076 & 9678897252 & 27130330208 & 22060899814 & 2169056374\\
20 & 814302292 & 23925376862 & 250129250080 & 1049599926148 & 1546833262014 & 510434205417\\
\midrule
\end{tabular}
\caption{Непомеченные $3$-регулярные карты}
\end{center}
\end{table}

\begin{table}[h!]
\begin{center}
\scriptsize
\begin{tabular}{c|cccccc}
\midrule
$ v \backslash g$ &\phantom{0000000}0\phantom{0000000}&\phantom{0000000}1\phantom{0000000}&\phantom{0000000}2\phantom{0000000}&\phantom{0000000}3\phantom{0000000}&\phantom{0000000}4\phantom{0000000}&\phantom{0000000}5\phantom{0000000}\\ 
\midrule
1 & 2 & 1 & 0 & 0 & 0 & 0\\
2 & 9 & 15 & 0 & 0 & 0 & 0\\
3 & 54 & 198 & 45 & 0 & 0 & 0\\
4 & 378 & 2511 & 2007 & 0 & 0 & 0\\
5 & 2916 & 31266 & 56646 & 9450 & 0 & 0\\
6 & 24057 & 385398 & 1290087 & 750762 & 0 & 0\\
7 & 208494 & 4721004 & 25872210 & 34001964 & 4729725 & 0\\
8 & 1876446 & 57590271 & 476712054 & 1155686130 & 555627195 & 0\\
9 & 17399772 & 700465482 & 8266938732 & 32793653844 & 35407720998 & 4341887550\\
10 & 165297834 & 8501284530 & 136971261063 & 820752362820 & 1628891511507 & 685100895750\\

\midrule
\end{tabular}
\caption{Помеченные $4$-регулярные карты}
\end{center}
\end{table}

\begin{table}[h!]
\begin{center}
\scriptsize
\begin{tabular}{c|cccccc}
\midrule
$ v \backslash g$ &\phantom{0000000}0\phantom{0000000}&\phantom{0000000}1\phantom{0000000}&\phantom{0000000}2\phantom{0000000}&\phantom{0000000}3\phantom{0000000}&\phantom{0000000}4\phantom{0000000}&\phantom{0000000}5\phantom{0000000}\\ 
\midrule
1 & 1 & 1 & 0 & 0 & 0 & 0\\
2 & 3 & 4 & 0 & 0 & 0 & 0\\
3 & 7 & 23 & 6 & 0 & 0 & 0\\
4 & 33 & 185 & 147 & 0 & 0 & 0\\
5 & 156 & 1647 & 2937 & 510 & 0 & 0\\
6 & 1070 & 16455 & 54511 & 31765 & 0 & 0\\
7 & 7515 & 169734 & 927117 & 1218026 & 169772 & 0\\
8 & 59151 & 1805028 & 14916436 & 36144608 & 17378795 & 0\\
9 & 483925 & 19472757 & 229711384 & 911103580 & 983694123 & 120644422\\
10 & 4136964 & 212603589 & 3424691109 & 20519920506 & 40723722416 & 17128196414\\
\midrule
\end{tabular}
\caption{Непомеченные $4$-регулярные карты}
\end{center}
\end{table}

\begin{table}[h!]
\begin{center}
\scriptsize
\begin{tabular}{c|cccccc}
\midrule
$ v \backslash g$ &\phantom{0000000}0\phantom{0000000}&\phantom{0000000}1\phantom{0000000}&\phantom{0000000}2\phantom{0000000}&\phantom{0000000}3\phantom{0000000}&\phantom{0000000}4\phantom{0000000}&\phantom{0000000}5\phantom{0000000}\\ 
\midrule
2 & 36 & 120 & 33 & 0 & 0 & 0\\
4 & 8640 & 125280 & 442800 & 292680 & 0 & 0\\
6 & 3312576 & 120800160 & 1457790048 & 6279955920 & 7582234716 & 1038647610\\
8 & 1586304000 & 113579366400 & 3119723797248 & 38667839616000 & 206853666065856 & 395867774517120\\
\midrule
\end{tabular}
\caption{Помеченные $5$-регулярные карты}
\end{center}
\end{table}

\begin{table}[h!]
\begin{center}
\scriptsize
\begin{tabular}{c|cccccc}
\midrule
$ v \backslash g$ &\phantom{0000000}0\phantom{0000000}&\phantom{0000000}1\phantom{0000000}&\phantom{0000000}2\phantom{0000000}&\phantom{0000000}3\phantom{0000000}&\phantom{0000000}4\phantom{0000000}&\phantom{0000000}5\phantom{0000000}\\ 
\midrule
2 & 7 & 15 & 7 & 0 & 0 & 0\\
4 & 468 & 6423 & 22409 & 14806 & 0 & 0\\
6 & 111096 & 4031952 & 48613214 & 209370256 & 252775238 & 34629024\\
8 & 39670362 & 2839677570 & 77994238389 & 966700391724 & 5171350497756 & 9896704783988\\
\midrule
\end{tabular}
\caption{Непомеченные $5$-регулярные карты}
\end{center}
\end{table}

\begin{table}[h!]
\begin{center}
\scriptsize
\begin{tabular}{c|cccccc}
\midrule
$ v \backslash g$ &\phantom{0000000}0\phantom{0000000}&\phantom{0000000}1\phantom{0000000}&\phantom{0000000}2\phantom{0000000}&\phantom{0000000}3\phantom{0000000}&\phantom{0000000}4\phantom{0000000}&\phantom{0000000}5\phantom{0000000}\\ 
\midrule
1 & 5 & 10 & 0 & 0 & 0 & 0\\
2 & 100 & 800 & 795 & 0 & 0 & 0\\
3 & 3000 & 58000 & 240150 & 171050 & 0 & 0\\
4 & 110000 & 4080000 & 41679000 & 124277000 & 72037775 & 0\\
5 & 4550000 & 283100000 & 5621625000 & 41509275000 & 100405619375 & 50323507500\\
6 & 204000000 & 19496000000 & 655512600000 & 9438042300000 & 56500699920000 & 117532756935000\\
\midrule
\end{tabular}
\caption{Помеченные $6$-регулярные карты}
\end{center}
\end{table}

\begin{table}[h!]
\begin{center}
\scriptsize
\begin{tabular}{c|cccccc}
\midrule
$ v \backslash g$ &\phantom{0000000}0\phantom{0000000}&\phantom{0000000}1\phantom{0000000}&\phantom{0000000}2\phantom{0000000}&\phantom{0000000}3\phantom{0000000}&\phantom{0000000}4\phantom{0000000}&\phantom{0000000}5\phantom{0000000}\\ 
\midrule
1 & 2 & 3 & 0 & 0 & 0 & 0\\
2 & 14 & 81 & 79 & 0 & 0 & 0\\
3 & 180 & 3313 & 13496 & 9625 & 0 & 0\\
4 & 4740 & 171282 & 1740813 & 5184959 & 3005225 & 0\\
5 & 152024 & 9444158 & 187427352 & 1383758744 & 3347008068 & 1677530177\\
6 & 5672485 & 541659909 & 18209426105 & 262170917920 & 1569470797305 & 3264807228317\\
\midrule
\end{tabular}
\caption{Непомеченные $6$-регулярные карты}
\label{table_last}
\end{center}
\end{table}

\end{document}